# MODEL BUILDING WITH MULTIPLE DEPENDENT VARIABLES AND CONSTRAINTS

Dr. Chris Tofallis
The Business School
University of Hertfordshire
U.K.

E-mail: c.tofallis@herts.ac.uk
Tel: +44 (0) 1707 285486

**Summary**

The most widely used method for finding relationships between several quantities is multiple regression. This however is restricted to a single dependent variable. We present a more general method which allows models to be constructed with *multiple* variables on both sides of an equation and which can be computed easily using a spreadsheet program. The underlying principle (originating from canonical correlation analysis) is that of maximising the correlation between the two sides of the model equation.

This paper presents a fitting procedure which makes it possible to force the estimated model to satisfy constraint conditions which it is required to possess, these may arise from theory, prior knowledge or be intuitively obvious. We also show that the least squares approach to the problem is inadequate as it produces models which are not scale invariant.

## 1. Introduction

Regression is one of the most widely used quantitative techniques in the natural and social sciences. Multiple regression allows one to fit a model to data when a single variable (y) is believed to depend on a number of other variables. Single equation models involving *multiple* dependent as well as multiple independent variables are much rarer in the literature. These might attempt to describe a relationship between a set of response variables and a set of explanatory variables, or between a set of outputs and a set of inputs. Here are some possible applications:

(i) We may be attempting to relate two quantities, which are not easily, or directly measurable by single variables so that some kind of index has to be constructed from several variables. Our approach permits an optimal set of weights to be found for combining the components of each index so that the strength of association between them is maximised. It is also possible to impose conditions on the model to force it to have certain desirable properties, or to ensure it abides by any restrictions it must adhere to. A simple example would be when we require particular weights to be positive and others negative, or when we require one side of the equation to always have a value in a certain range e.g. 0 to 100%.

(ii) Consider a system for which there are a number of variables under the control of the experimenter and also a number of other response variables, which are measured. In each run of the experiment a different set of values is chosen for the controlled variables, and the response variables are recorded. The proposed method could be used on the data collected to try to find a model, which relates all the variables. (Of course the situation can also be one which is being observed rather than an experimental design e.g. levels of various pollutants in a river).

(iii) Consider a set of objects or animals or people, which have some of their attributes (e.g. physical characteristics), measured. Measurements are also taken of some of their abilities or functions. The method could be used to relate these two sets of quantities.

(iv)Another possible application is in discovering resonances (also called commensurabilities) in a set of observations e.g. from the physical sciences or engineering. For instance, in celestial mechanics the shape and position of an orbit is naturally described by a set of variables known as orbital elements. Very often two orbiting bodies (e.g. Jupiter and an asteroid) are found to exist in a resonance i.e. with a linear combination of the elements of one orbit equating to a linear combination of the elements of the second orbit, with the coefficients being integer. The method of this paper could be applied to find the combinations which give maximum correlation. If these contained coefficients which were close to small integers then a resonance may be detected. Large integers are not of interest as they are not connected with strong repeated forcing effects. One can impose upper limit constraints on the size of these coefficients prior to the model being constructed. An alternative approach would be to find the maximum correlation combinations subject to the condition that only integer coefficients be permitted. This integer optimisation can also be solved on a spreadsheet although such problems (integer programming problems) are known to be computationally more time consuming.

In general the attraction of a single equation model is that it may provide a compact way of bringing together all the relevant factors. This paper describes how such a model can be very simply constructed using nothing more than a spreadsheet package. As the breadth of application of this technique is envisioned to be very wide we have deliberately used very simple language to permit researchers in all fields to follow the procedure.

## 2. Canonical correlation analysis

Among the methods one finds in textbooks on multivariate statistics (e.g. Hair et al 1995, or Tabachnick and Fidell 1996) there is a surprisingly little used technique known as canonical correlation analysis (CCA). This finds a linear combination (i.e. a weighted sum) of the dependent variables and a linear combination of the independent variables such that the correlation between the two is maximised. According to Tabachnick and Fidell (1996, p.196) CCA is 'one of the most general of multivariate techniques...but it is also the least used and most impoverished'. One reason CCA has not been widely used is that the values of the weights (coefficients, or parameter values) which arise are often intuitively unacceptable e.g. they may have the wrong sign or be of the wrong relative magnitude. This author was unable to find any statistical software package offering CCA which could avoid these difficulties by the imposition of constraints, moreover the problem appears not to have been tackled in statistical journals. This is the motivation for the present paper.

In order to compute the weights/coefficients CCA solves an eigenvector problem. The procedure is described in Tabachnick and Fidell (1996) where the authors concede that 'it is not particularly enlightening' (p.201). Furthermore it is not apparent how it would be possible to impose conditions on the weights *a priori* in such an eigenvector approach. We circumvent this difficulty by re-casting the problem as a constrained optimisation i.e. find values for the weights so as to maximise the correlation between the two linear combinations subject to any conditions one needs to impose.

**3. Constrained optimisation on a spreadsheet**

The most widely used spreadsheet packages now have a built-in facility or 'tool' for constrained optimisation, this is often called a 'solver' or 'optimiser'. The author uses the Microsoft Excel spreadsheet program but the providers of its solver tool (Frontline Systems) also provide it on the Lotus 123 and Quattro Pro programs. We shall explain how to use this tool for our purposes, but first we describe how the spreadsheet may be set up in preparation.

Assuming the data appears on the spreadsheet with one variable per column, we need to allocate a row of cells to hold the associated weights, say at the bottom of the data columns (see Figure 1). We also set up two columns to hold the values of the weighted combinations (call them X and Y).

Figure 1. An example showing the spreadsheet layout when there are two x-variables and three y-variables. The values would appear below their indicated labels.

| $x_1$ | $x_2$ | $X = a_1x_1 + a_2x_2$ | $y_1$ | $y_2$ | $y_3$ | $Y = \Sigma by$ |
|---|---|---|---|---|---|---|
| | | | | | | |
| $a_1$ | $a_2$ | | $b_1$ | $b_2$ | $b_3$ | |
| | | | | | | |

Finally we shall need a cell (outside the table) which calculates the correlation between X and Y. The correlation function is available as a standard function in spreadsheets. It is always convenient and more meaningful to use the spreadsheet's facility to attach names to particular cells or groups of cells since one can then refer to them by name rather than by cell address. For instance it will be useful to attach the name 'Weights' to the row of cells we have set aside for them (the bottom row).

Each constraint (be it an inequality or an equation) needs to be set up by placing an expression for evaluating its left hand side in one cell and the right hand side in another cell (not shown in Figure 1). Traditionally the constraints are arranged so that all terms containing variables are on the left, so that the right hand side merely contains a constant, which may be zero. The solver tool can now be selected. Within its dialog box you will be asked to specify the 'target cell' or objective function - this is the cell containing the quantity to be maximised i.e. the correlation. Next you need to enter the 'changing cells' or decision variables; these will contain the values of the weights and can be specified most simply by typing in the name 'Weights' which we have already defined. Lastly we input the constraints. They are entered by giving the cell addresses for the left and right hand sides of the constraints and specifying their type ($\geq$, =, or $\leq$) e.g. for non-negativity we simply enter 'Weights $\geq$ 0'. The 'options' dialogue box permits various settings to be made. The most important are firstly not to select 'assume a linear problem', and secondly to select 'automatic scaling', the latter reduces round-off error when there are large differences in the magnitude of data values. The most recent version of Excel (1997) allows one to set a value for a 'convergence' parameter. This is expressed as the relative change (between 0 and 1) in the quantity being optimised, if it does not change by this amount over five iterations then the Solver stops and returns its results. (The actual optimisation method used is the generalised reduced gradient procedure.) When the solving process is finished it is important to note the completion message that Solver provides. If it does not say that the optimality conditions are satisfied then one should reduce the convergence parameter value and re-solve until optimality is confirmed, this of course requires a longer computation time. (A discussion relating to global optimality appears in the appendix.) The optimal values for the weights and the maximum correlation will appear in the cells that we have set aside for them.

The correlation between two quantities is unaffected if either one or both are multiplied by a constant. This implies that we can do this to one or both of our sets of weights. For instance it may be convenient to scale the weights so that one of them takes the value of unity. In fact such a normalisation condition can be used as a constraint in the above solution procedure in order to avoid the trivial solution where all the weights are assigned a value of zero. An alternative normalisation might be that one set of weights should sum to unity. There is no loss of generality arising from this.

In order to check results generated by the spreadsheet with those from specialist statistical software which handles canonical correlation analysis we used the worked example from Tabachnick and Fidell (1996), which of course does not involve any constraints. The spreadsheet provided equivalent results, to within a scaling factor, to results from the programs SPSS, SYSTAT, SAS, and BMDP6M. It is worth noting that according to the manual the solver in Version 5 of Excel is able to deal with 100 constraints and 200 parameters.

**4. Building the model**

If the value of the correlation that has been found is high then it follows that a straight-line relationship between the composite variables X and Y will provide a good fit. So far the procedure we have followed has treated the dependent and independent variables in the same way but once we proceed to find a regression equation then this symmetry will be lost: a regression of Y on X will not be equivalent to a regression of X on Y. (Conventional regression minimises the sum of the squared vertical (Y) deviations from the fitted model. One way of avoiding the asymmetry might involve taking the perpendicular distance to the graph instead, see Gander and von Matt (1993) for details).

If however the correlation value between X and Y is too small then to improve matters the model can be enriched by adding variables which are functions of the existing ones e.g. squared terms or logarithms to help deal with nonlinearities, or products of two or more of the original variables to help deal with interactions between factors. These additional quantities will appear as further columns of data and will be treated in exactly the same way as the original data - because they are still multiplied by a simple weight there is no additional complexity involved.

The above procedure is then repeated to find the new sets of coefficients *a* and *b* which maximise the correlation between Σ a f(x) and Σ b g(y).

## 5. Comparison with least squares

How do the models generated by our approach compare with least squares applied in the following way?

$$\left(\Sigma\, b\, y - \Sigma\, a\, x - c\right)_i = e_i \qquad (1)$$

(where *i* relates to a particular row in the data table, see Figure 1), and find values of the coefficients to minimise $\Sigma e_i^2$.

Once again a normalisation condition is required in order to avoid obtaining a solution in which all the coefficients are zero, giving $\Sigma e_i^2 = 0$. One possibility is to set one of the coefficients to unity, this merely implies multiplying through the model equation by a suitable constant. Another way of looking at this is that appropriate units of measurement are chosen so that the selected coefficient is unity. This might at first seem a harmless thing to do, after all we would not expect a model to be affected by the units of measurement, rather we would expect it to be 'scale invariant' or 'units invariant'. We shall see that this is not the case.

The results of a comparison between the two approaches are most illuminating. Firstly, the least squares approach does not provide a model with maximum correlation between Y and X, the correlation is generally lower and, by definition, cannot be higher. Secondly, we get different models according to which coefficient has been normalised; these models are *not* equivalent. This can be understood as follows: setting say $b_k = 1$, means that (1) can be re-written as:

$$y_k = \Sigma\, a\, x + c - \sum_{j \neq k} b\, y + e$$

This is precisely least squares multiple regression with $y_k$ as the dependent variable. Now suppose instead that we set a different coefficient to unity. We are then changing the dependent variable. It is well known that changing the dependent variable gives rise to a different, non-equivalent model. This is true even with simple bivariate regression (one x-variable, one y-variable): regressing y on x is not the same as regressing x on y. What multiple regression does achieve is to maximise the correlation between the observed and predicted values of the dependent variable (the latter being a linear combination of the remaining variables); this is clearly not the same as maximising the correlation between X and Y.

In summary the approach of this paper, namely *maximum correlation modelling*, is superior to the least squares approach described above in that it provides equivalent models irrespective of the units of measurement of the individual variables used. Maximum correlation modelling thus has the useful property of being scale invariant i.e. changing the units of measurement of any variable does not affect the final model (the associated coefficient will simply scale accordingly to generate an equivalent model).

## 6. Example application

The data we shall use has been adapted from Ganley and Cubbin (1992, chapter 3). The units of analysis are 96 English Local Education Authorities (LEAs). The aim is to relate the examination results from secondary (high) school pupils in schools maintained by these authorities to a number of explanatory variables. The latter are sometimes called contextual or environmental variables. The examination results we use are the averages over a period of three years:1980/81 to 1982/3. At that time there were two types of examination in use for sixteen-year-olds: ordinary ('O') levels and the lower level CSE (Certificate of Secondary Education), with the top grade of the latter being equivalent to an O level.

Three outcome variables will be used:

$y_1$: The percentage of pupils achieving at least 5 higher grade passes at O level/CSE,

$y_2$: The percentage of pupils achieving at least 6 graded results at O level /CSE, excluding pupils included in $y_1$,

$y_3$: The percentage of pupils achieving at least one graded result, excluding pupils already counted in $y_2$ and $y_1$.

The mean values of these outcomes across the 96 LEAs were 22.6%, 40.5%, and 25.5% respectively. The mean percentage across the LEAs achieving no passes at all was 11.3%.

The explanatory variables used are:

$x_1$: Secondary school teaching expenditure per pupil.

$x_2$: Percentage of pupils living in households whose head is a non-manual worker, excluding junior non-manual workers and non-manual supervisors.

$x_3$: Percentage of pupils living in households, which have all the standard amenities and are not overcrowded (density less than 1.5 persons per room).

$x_4$: Percentage of pupils born in the UK, Ireland, USA, the old Commonwealth, or whose head of household was born in any of these places.

$x_5$: Persons per hectare.

$x_6$: (Persons per hectare)$^2$. This quadratic term was included because a monotonic relationship between exam performance and population density is unlikely; whilst the high densities associated with inner cities might produce a negative effect, so might the very low densities in rural regions, this may arise due to lack of facilities or resources.

We seek some measure of exam performance which is a weighted combination of the three outcome variables: $Y = b_1 y_1 + b_2 y_2 + b_3 y_3$. It would be preferable if the weights in some way reflected the fact that these are three distinct levels of achievement. It is clearly more difficult for a pupil to achieve at least five higher grade passes and so be included in $y_1$, than to be counted under $y_2$ or $y_3$. Similarly it is easier to be counted in the $y_3$ group than in $y_2$. Hence we shall impose the following inequality constraints on the weights: $b_1 \geq b_2 \geq b_3$. Some spreadsheet solvers require that such conditions be entered with a constant on the right hand side, so we have $b_1 - b_2 \geq 0$ and $b_2 - b_3 \geq 0$. No conditions were imposed on the x coefficients.

We then solve to find linear combinations $X = \Sigma ax$ and $Y = \Sigma by$ which have maximum correlation whilst satisfying the above constraints. Using a Pentium processor running at 233megahertz the computation time was less than five seconds. In the results which follow we have scaled the coefficients in Y to make $b_2 = 1$ for convenience, there is no loss of generality in doing this):

$$Y = 2.871 \, y_1 + y_2 + y_3 \tag{2}$$

$$X = 0.0071\, x_1 + 0.471\, x_2 + 0.432\, x_3 - 0.0083\, x_4 + 0.1007 x_5 - 0.0025\, x_6$$

The correlation between X and Y is 0.9023 .

There are some interesting observations one can make from these results. Firstly, the small negative coefficient on $x_4$ indicates that a high proportion of pupils from *outside* the UK, Ireland, USA and the old Commonwealth actually acts to improve expected performance. Other studies have shown a similar effect. For instance Nuttall et al (1989) carried out a detailed study which included the categorisation of pupils into 11 different ethnic backgrounds; they found that all but one of the "ethnic groups perform significantly better than the ESWIs" (pupils of English, Scottish, Welsh or Irish background). Secondly, when the quadratic function of population density ($a_5 x_5 + a_6 x_6$ or $a_5 x_5 + a_6 x_5^2$ ) is plotted it shows a peak at 40 persons per hectare; the further the density is from this level (either above or below) - the lower the expected score. When the analysis was repeated without any constraints i.e. a straightforward canonical correlation analysis, it was found that the coefficient of $y_3$ exceeded that of $y_2$, which most would consider unacceptable in any index of performance. Removing the constraints will generally improve the correlation value. In this case it was a negligible change, from 0.9023 to 0.9024. Thus including the constraints allowed us to generate a more intuitively acceptable model for next to no loss in the goodness of fit. However our constrained model is not entirely satisfactory in that two of the outcome coefficients are equal. If this were considered unsatisfactory one could impose an additional constraint that forced a minimal acceptable difference between these two weights. The difficulty with this is in deciding what the minimal difference should be. An alternative approach would be to specify an acceptable value for the correlation and then deduce the weights subject to the existing constraints. This is easily achieved using a facility within the solver which allows one to set a value for the target cell rather than to maximise it. When a target correlation of 0.900 was set the resulting outcome weights were in the correct order: 2.84, 1, and 0.886 respectively, this was at a cost of only about two parts in 900 in the correlation. The weights on the x variables were hardly affected by this change.

A regression of Y on X using (2) produced: Y = 2.35 X −7.01 (with $r^2$ = 0.81). Such a formula, when expanded out, relates in a compact form all the variables we have been dealing with. It could for instance be used to estimate the expected performance in a given LEA for the given levels of the environmental variables. This could then be compared with actual performance.

We must stress that the application we have described above is for illustrative purposes only. One can criticise it on a number of grounds, from the choice of variable used to the form of the model. No doubt it can be improved in a number of ways.

Nowhere have we commented on the statistical significance of the coefficients or the correlation value obtained from the method we have presented. However it should be noted that we have used a complete data set not merely a sample of LEAs. In fact the necessary inference theory remains to be constructed. We have presented a fitting procedure alone, we are not aware of any work, which allows statistical inference to be carried out for constrained canonical correlation analysis.

## 7. Conclusion

When modelling relationships between quantities which fall into two groups or classes (e.g. inputs and outputs, or environmental and response, or independent and dependent), the simplest model would be a single equation with the variables belonging to each group appearing on separate sides of the equation. We have presented a simple approach for fitting such a model based on finding the parameter values or weights, which maximise the correlation between the two sides of such an equation. In addition the procedure allows us to include any *a priori* information we have (e.g. relationships which are known to apply) so that the model will have the required properties. We have also shown that a least squares approach to such a fitting problem generates models which are not scale invariant, whilst the proposed procedure does not suffer from this inadequacy. No specialist software program is required as the procedure can be carried out on a spreadsheet. Researchers will thus be able to make immediate use of this tool (which we prefer to call '*maximum correlation modelling*' rather than the more cumbersome 'constrained canonical correlation analysis'

A future paper will show how this tool can be applied to the problem of setting up a measure of value-added in an educational setting, whilst in Tofallis (1997, 1998) the method is combined with data envelopment analysis to show how best-practice can be modelled.

**Appendix**

In this appendix we discuss under what circumstances the solution found by the optimisation process is a global optimum rather than merely a local optimum. In our context a local optimum that is not a global optimum would mean that a higher value of the correlation could be found if the optimisation process searched in a different part of the space of possible solutions (also known as the feasible region). A similar difficulty arises in nonlinear regression (where the fitted parameters do not appear as simple linear coefficients). The traditional advice there is to start the search process using values of the parameters which experience or theory suggest are reasonable, and once a solution is found to experiment with other starting values as a precautionary measure. We shall show that this will not be necessary i.e. that our optimum will be global, provided that the set of constraint conditions enclose a convex region (this means that a line joining any two points in the region will lie entirely within the region). If all the constraints are linear functions of the weights (non-negativity conditions on the weights are the simplest example) then the region is guaranteed to be convex. Of course it is also possible to have a convex region where the boundaries are not flat planes. These conditions will apply in the majority of practical applications.

Our demonstration will make use of the following basic result from optimisation theory (Taha, 1992): minimising a convex function over a convex region will always give a global optimum, (a convex function is one where a line joining two points on the function surface or graph will never pass below that surface). Since we are assuming that we are dealing with a convex region we need to show that our problem is equivalent to minimising a convex function.

If we denote the correlation between X and Y by r(X,Y) then we have that r(X,Y) = −r(X, −Y). It follows that maximising r(X,Y) is equivalent to minimising r(X, −Y).

Let $X = \Sigma\ u_i\ x_i$ and $Y = \Sigma\ v_j\ y_j$ where the $u_i$ and $v_j$ are the weights attached to the independent and dependent data variables respectively. From the definition of correlation it follows that r(X, −Y) is a quadratic form in the weight variables. Let us assume that:

(i) We have arranged that the weights will all be positive; this can be done by multiplying particular variables by −1 if necessary,

(ii) An appropriate constant has been added to each variable to make all its values positive. (These last two steps constitute linear transformations of the data and so there is no loss of generality - one can always go back and reverse the transformations to recover the original variables when the optimal model has been found.) These transformations ensure that r(X, −Y) is negative definite. A negative definite quadratic form is a convex function (Taha. 1992) and so our optimisation will provide the required global optimum.

**Acknowledgement**
I am grateful to R. Ian Cooper for introducing me to canonical correlation analysis.